\documentclass[12pt]{article}
\usepackage{graphicx}
\usepackage{color}
\usepackage{amsmath,amsthm}
\usepackage{amssymb}
\usepackage{xspace}
\usepackage{epsfig}
\usepackage[dvipsnames]{xcolor}
\usepackage{soul}
\usepackage{xfrac}
\usepackage{hyperref}

\usepackage{enumerate}

\usepackage{authblk}

\usepackage[margin=1in]{geometry}

\newtheorem{theorem}{Theorem}[section]
\newtheorem*{maintheorem}{Theorem \ref{thm:level_2_flip-connectivity}}
\newtheorem{lemma}[theorem]{Lemma}
\newtheorem{corollary}[theorem]{Corollary}

\newtheorem{proposition}[theorem]{Proposition}
\theoremstyle{definition}
\newtheorem{definition}[theorem]{Definition}
\newtheorem{example}[theorem]{Example}

\theoremstyle{remark}
\newtheorem*{remark}{Remark}

\newcommand {\mm}[1] {\ifmmode{#1}\else{\mbox{\(#1\)}}\fi}

\newcommand {\ceiling}[1] {{\left\lceil #1 \right\rceil}}

\newcommand{\ignore}[1]{}

\newsavebox{\smallProofsym}                 
\savebox{\smallProofsym}                               %
{
\begin{picture}(6,6)
\put(0,0){\framebox(6,6){}}
\put(0,2){\framebox(4,4){}}
\end{picture} 
}  

\makeatletter
\long\def\@makecaption#1#2{%
  \vskip\abovecaptionskip
  \sbox\@tempboxa{\small #1: #2}%
  \ifdim \wd\@tempboxa >\hsize
    \small #1: #2\par
  \else
    \global \@minipagefalse
    \hb@xt@\hsize{\hfil\box\@tempboxa\hfil}%
  \fi
  \vskip\belowcaptionskip}
\makeatother

\theoremstyle{plain}

\newtheorem*{supermaintheorem*}{Main Theorem}

\newtheorem*{supermaincorollary*}{Main Corollary}

\DeclareMathOperator{\area}{area}

\newcommand{\Rspace}        {\mm{{\mathbb R}}}

\newcommand{\Level}[2]      {\mm{{#1}^{({#2})}}}

\newcommand{\GKZ}[1]        {\mm{{\rm GKZ}{({#1})}}}

\newcommand{\card}[1]       {\mm{{|}{#1}{|}}}
\newcommand{\conv}[1]       {\mm{{\rm conv}{#1}}}

\newcommand{\Skip}[1]       {}

\definecolor{blue-green}{rgb}{0.0, 0.87, 0.87}

\begin{document}

\title{Flips in Two-dimensional Hypertriangulations\footnote{
{\bf MSC classes:} 52C20, 52C40.}}

\author{Herbert Edelsbrunner\footnote{ISTA (Institute of Science and Technology Austria), Kloster\-neu\-burg, Austria, \texttt{edels@ist.ac.at}}, Alexey Garber\footnote{School of Mathematical and Statistical Sciences, University of Texas Rio Grande Valley, Brownsville, Texas, USA, \texttt{alexey.garber@utrgv.edu}}, Mohadese Ghafari\footnote{Khoury College of Computer Science, Northeastern University, Boston, Massachusetts, USA, \texttt{ghafari.m@northeastern.edu}}, 
\\Teresa Heiss\footnote{The Australian National University, Canberra ACT, Australia, \texttt{teresa.heiss.synak@anu.edu.au}}, Morteza Saghafian\footnote{ISTA (Institute of Science and Technology Austria), Kloster\-neu\-burg, Austria, \texttt{morteza.saghafian @ist.ac.at}}}

\maketitle

\begin{abstract}
  We study flips in hypertriangulations of planar points sets.
  Here a level-$k$ hypertriangulation of $n$ points in the plane is a subdivision induced by the projection of a $k$-hypersimplex, which is the convex hull of the barycenters of the $(k-1)$-dimensional faces of the standard $(n-1)$-simplex.
  In particular, we introduce four types of flips and prove that the level-2 hypertriangulations are connected by these flips. 
\end{abstract}

\section{Introduction}
\label{sec:1}

Triangulations appear in many fields of pure and applied mathematics, and they go back to the study of the Catalan numbers by Euler and Goldbach in the early 18th century \cite{Pak14}, if not further.
Flips were introduced by Wagner \cite{Wag36} as a tool to study how triangulations change.
In particular, he proved that for a finite planar set, the family of triangulations is flip-connected.
This fact was later exploited by Lawson \cite{Law77} to construct triangulations algorithmically for the purpose of surface interpolation.
Ever since, flip-connectivity has become a standard topic in discrete and computational geometry.
In the plane, flip-connectivity leads to efficient algorithms for constructing Delaunay triangulations and to proofs of extremal properties, for example that among all triangulations of a finite set the Delaunay triangulation maximizes the minimum angle \cite{Sib78}.

\medskip
Beyond two dimensions, the connectivity with flips---also called \emph{Pachner moves} or \emph{bistellar flips}---is a more challenging concept.
For example, the greedy algorithm that flips a locally non-convex configuration succeeds in constructing the Delaunay triangulation in the plane \cite{Law77}, but can get stuck in local minima in three dimensions; see \cite{Joe89} for examples.
However, a more limited protocol that inserts points incrementally and repairs the Delaunay triangulation after every insertion also succeeds in three dimensions \cite{Joe91} and extends to higher dimensions and to coherent triangulations \cite{EdSh96}.
Note that the latter are known in the literature under a variety of names, including \emph{Laguerre}, \emph{regular}, and \emph{weighted Delaunay triangulations}.
While flip-connectivity in three dimensions is still an open question, Santos \cite{San00} has exhibited a configuration in six dimensions whose family of triangulations is not flip-connected.
We refer to \cite{EdRe97,HNU99, Liu18, Liu20, San05, WaWe22} and references therein
for a multitude of results on flip-connectivity in a variety of settings.

\medskip
In a more general setting beyond triangulations, flips were introduced in the study of subdivisions induced by projections of polytopes and associated fiber polytopes by Billera and Sturmfels \cite{BiSt92}. 
The corresponding notion of flips comes from the Baues poset, see \cite{BKS94,Rei99,San01}. It seems that such flips were first discussed by Billera, Kapranov, and Sturmfels in \cite{BKS94} as possible generalizations of bistellar moves in the setting of the generalized Baues conjecture and related connectivity.
We refer to these as \emph{combinatorial} or \emph{Baues flips} as it seems difficult to give a geometric description in all possible cases, even for a generic projection. 
Primarily in connection to the Baues problem, the connectivity of the flip graph for various settings related to the subdivisions induced by projections of polytopes has been studied in \cite{Ath01, BKS94, Liu20a, RaZi96}, where positive as well as negative results are described.

\medskip
In this paper, we study flips in the family of hypertriangulations of a finite set of $n$ points in the plane. The term {\it hypertriangulations}, or more generally, {\it hypersimplicial subdivisions}, was introduced by Olarte and Santos~\cite{OlSa21} to study the Baues problem for families of plabic graphs related to the totally positive Grassmannian \cite{Pos18}; these constructions and their connections to plabic graphs appeared in earlier works by Postnikov \cite{Pos18} and Galashin \cite{Gal18}. 
Hypertriangulations are triangulations induced by projections of an $(n-1)$-dimensional hypersimplex to the plane.
To explain these concepts, we fix an integer $k$ between $1$ and $n-1$, called the \emph{level}, and we write $\Delta_n^{(k)}$ for the $k$-fold scaling of the convex hull of the barycenters of the $(k-1)$-dimensional faces of the standard $(n-1)$-simplex, $\Delta_n = \Delta_n^{(1)}$.
Correspondingly, we write $\Level{A}{k}$ for the set of $k$-fold sums of the points in a given set of $n$ points, $A = \Level{A}{1}$.
The projection fixes a bijection between the vertices of $\Delta_n$ and the points of $A$ and, by construction, maps the vertices of $\Delta_n^{(k)}$ to the points of $\Level{A}{k}$.
The hypertriangulations follow by selecting and projecting appropriate subsets of the $2$-dimensional faces of $\Delta_n^{(k)}$.
The subclass of coherent level-$k$ hypertriangulations are also known as \emph{weighted order-$k$ Delaunay triangulations}, 
defined by generalizing order-$k$ Delaunay triangulations \cite{Aur90} to the weighted setting as suggested in \cite{EdOs22}\footnote{More formally, given a set of points in $\Rspace^2$ with assigned real weights, we lift the points to paraboloid $y=x_1^2+x_2^2$ in $\Rspace^3$ and shift each point vertically according to its weight. Considering the convex hull of all $k$-fold sums of the lifted points we obtain the order-$k$ Delaunay triangulation of the original weigthed point set. Varying the weights, we may obtain all coherent hypertriangulations as described in Section \ref{sec:2.3} below. An alternative way to obtain order-$k$ Delaunay triangulations is using duals to $k$th-order power diagrams (weighted order-$k$ Voronoi diagrams), see \cite[Sect. 13.6]{Ede87}.}.
The case $k=1$ corresponds to the family of usual triangulations of the set $A$ that can be viewed as (tight) subdivisions induced by the projection of (the vertices of) the standard simplex $\Delta_n$ on the points of $A$.
Our main result establishes flip-connectivity of level-2 hypertriangulations of generic point sets, as stated and proved in Section~\ref{sec:4}.

\medskip
\begin{maintheorem}
  For every generic point set $A \subseteq \Rspace^2$, the level-$2$ hypertriangulations of $A$ are flip-connected.
\end{maintheorem}

\medskip
We note that when $A$ is in convex position, then flip-connectivity of the family of all level-$k$ hypertriangulations of $A$ can be deduced from results of Postnikov \cite[Cor. 11.2]{Pos18}.

\medskip
The flips we consider in the theorem above are possibly more restrictive than the Baues flips. They are introduced geometrically in Section \ref{sec:3}. For the specific setting of projections of hypersimplices, these flips coincide with \emph{$\pi$-flips} introduced  by Santos in \cite{San01} 
and are subsets of the Baues flips for the corresponding projections. However, we were unable to prove or disprove that the family of flips we use coincides with the family of Baues flips for (generic) projections of hypersimplices.

\medskip
Since all two-dimensional faces of hypersimplices are triangles, the family of level-$k$ hypertriangulations of $A$ is a subset of the family of triangulations of the set $A^{(k)}$. Notwithstanding, the flip connectivity of the former family does not follow from the flip-connectivity of the latter family and the flips are defined differently. Moreover, flip-connectivity using Baues flips is not guaranteed for projections of simplicial polytopes on the plane as illustrated by the example of Rambau and Ziegler \cite{RaZi96}. Thus, in this paper we prove that every (generic) projection of the hypersimplex $\Delta_n^{(2)}$ gives rise to a connected Baues poset extending the same property known for projections of simplices.

\medskip
Before presenting the outline of this paper, we stress that we restrict ourselves to generic point sets in the plane.
Section~\ref{sec:2} follows Olarte and Santos \cite{OlSa21} and introduces level-$k$ hypertriangulations of a set of $n$ points, $A$, as tight subdivisions induced by the projection of the hypersimplex $\Delta_n^{(k)}$. 
In particular, we give a combinatorial description of such triangulations without using the associated projection.
We also define the more general hypersimplicial subdivisions and sketch a connection to fiber polytopes \cite{BiSt92}.
Section~\ref{sec:3} introduces the four types of flips for level-$k$ hypertriangulations. 
Notably, these flips are geometric, they do not depend on the level $k$, and they include (colored versions of) the classic flips for triangulations of planar point sets \cite{DRS10}.
Section~\ref{sec:4}, and Section \ref{sec:4.3} in particular, proves the main result of this paper: that the family of level-2 hypertriangulations of every (generic) planar point set is flip-connected. 
The main tool in the proof is the aging function for triangles of hypertriangulations defined in Section~\ref{sec:4.1}. 
We note that the aging function was used before by Galashin \cite{Gal18} and Olarte and Santos \cite{OlSa21} for point sets in convex position; particularly, Galashin \cite[Cor. 4.4]{Gal18} (see also \cite[Lem. 3.6]{BaWe20}) proved that the aging function and its inverse are well-defined for all hypertriangulations and all suitable levels $k$ if $A$ is in convex position. 
For our approach, properties of the aging function for arbitrary generic point sets are instrumental and we prove these properties in Sections~\ref{sec:4.1} and~\ref{sec:4.2}. In Section~\ref{sec:4.4}, we adapt an example by Olarte and Santos \cite{OlSa21} to show that the aging function may fail in general. In Section~\ref{sec:4.5}, we briefly discuss coherent hypertriangulations and their flip-connectivity.
%
%
Section~\ref{sec:6} concludes the paper with a discussion of possible further research.

\section{Introduction to Hypertriangulations}
\label{sec:2}

This section explains the main object of study: the hypertriangulations of a finite point set.
To begin, we give an informal introduction to the subject, drawing connections to the theory of fiber polytopes along the way, and follow up with formal definitions thereafter.

\subsection{Level-$k$ Hypertriangulations}
\label{sec:2.1}

Let $A = \{a_1, a_2, \ldots, a_n\}$ be a set of $n$ points in $\Rspace^2$.
Write $[n] = \{1,2,\ldots,n\}$, and for a subset $I \subseteq [n]$, let $a_I = \sum_{i \in I} a_i$ be the vector sum of the points with index in $I$.
Fixing a parameter $k \in [n-1]$, we write $\Level{A}{k} = \{ a_I \mid I \subseteq [n], \card{I}=k \}$ for the $k$-fold sums and consider all \emph{partial triangulations} of $\Level{A}{k}$, by which we mean the decompositions of the convex hull of $\Level{A}{k}$ into triangles with vertices in $\Level{A}{k}$ such that no vertex of one triangle lies on a side of another triangle.
We assume that $A$ is {\it generic} and for the time being, this means that no three points of $\Level{A}{k}$ are collinear. We will weaken that restriction to a standard one in Section \ref{sec:2.2}.
\begin{definition}
  \label{def:level-k_hypertriangulation}
  A \emph{level-$k$ hypertriangulation of $A$} is a partial triangulation of $\Level{A}{k}$ so that
  \begin{itemize}
    \item[(V)]
      every vertex is of the form $a_I$, with $\card{I} = k$, and
    \item[(E)]
      every edge connects two vertices, $a_I$ and $a_J$, that satisfy $\card{I \cap J} = k-1$.
  \end{itemize}
\end{definition}

\begin{figure}[hbt]
  \centering
  \vspace{0.0in}
  \resizebox{!}{2.2in}{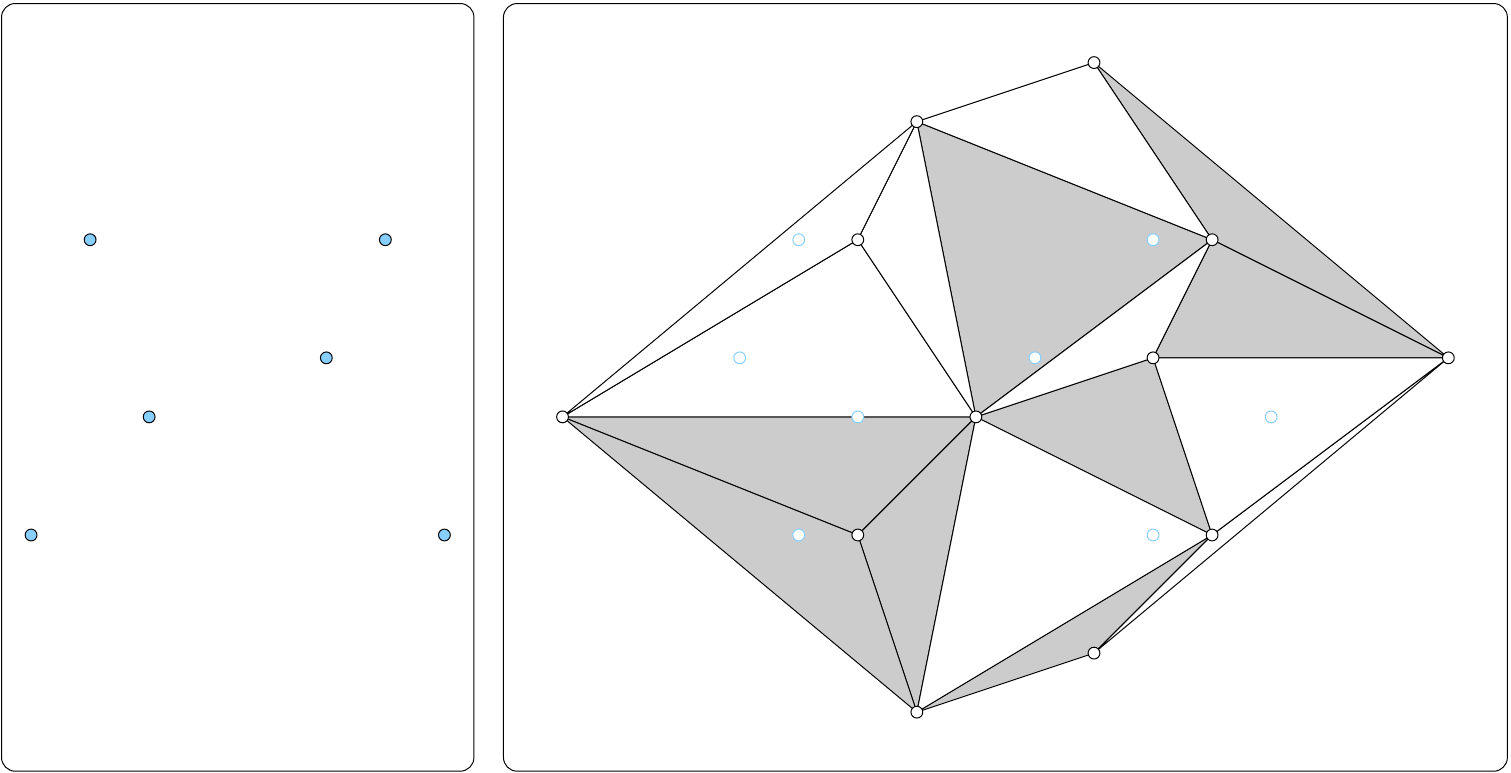}
    \caption{\emph{Left:} $n = 6$ points in the plane.
    \emph{Right:} the $\binom{n}{3} = 20$ triple sums of these points together with a triangulation of $12$ of the $20$ points.
    The three vertices of every \emph{black} triangle share one point in their sums, and the three vertices of every \emph{white} triangle share two.}
    \label{fig:alltypes}
\end{figure}
See Figure~\ref{fig:alltypes} for a level-$3$ hypertriangulation of $6$ points in $\Rspace^2$ as an example.
The labels of the vertices indicate how they are obtained from the $6$ points shown on the left.
Note that some of the $20$ points in $\Level{A}{3}$ are not vertices of the displayed hypertriangulation.
The requirement on the endpoints of every edge implies a similar requirement on the vertices of every triangle:
\begin{definition}
  \label{def:white_and_black_triangles}
  Let $\Delta=a_Ia_Ja_K$ be a triangle whose vertices and edges satisfy conditions (V) and (E).
  Then either $|I \cap J \cap K| = k-1$, in which case we call $\Delta$ a  \emph{white triangle}, or $|I \cap J \cap K| = k-2$, in which case we call $\Delta$ a \emph{black triangle}.
\end{definition}
Note that white triangles exist for $1 \leq k \leq n-2$ and black triangles exist for $2 \leq k \leq n-1$.
For a given triangulation, $T$, we write $W(T)$ and $B(T)$ for the sets of white and black triangles, respectively.
In the example in Figure~\ref{fig:alltypes}, there are $8$ triangles of each color.

\subsection{Hypersimplicial Subdivisions}
\label{sec:2.2}

The following interpretation of the above concepts gives an equivalent description within the theory of fiber polytopes and induced subdivisions.
We refer to \cite[Chap. 9]{Zie08} and \cite{DRS10} for a comprehensive introductions to this theory.

\medskip
Write $\Delta_n = \conv{\{e_1, e_2, \ldots, e_n\}} \subseteq \Rspace^n$ for the standard $(n-1)$-simplex,
and more generally $\Delta_n^{(k)} = \conv{\{e_I \mid I \subseteq [n], \card{I} = k\}}$ for the $k$-th standard $(n-1)$-dimensional hypersimplex.
Let $\pi \colon \Delta_n \to \conv{A}$ be the (linear) projection defined by $\pi (e_i) = a_i$, and following Olarte and Santos \cite{OlSa21}, we extend this to the projection $\pi \colon \Delta_n^{(k)} \to \conv{\Level{A}{k}}$.
We get tilings of $\Level{A}{k}$ by projecting subsets of the set of $2$-dimensional faces of $\Delta_n^{(k)}$ instead of the entire hypersimplex.
By construction, the label of each vertex of $\Delta_n^{(k)}$ is a subset of $k$ integers in $[n]$, the endpoints of each edge have labels that differ in exactly one integer, and the $2$-dimensional faces are triangles and therefore satisfy the condition on the vertex labels given in Definition~\ref{def:white_and_black_triangles}.
Each such tiling is called a \emph{hypertriangulation} in \cite{OlSa21}, and it agrees with our notion of hypertriangulation given in Definition~\ref{def:level-k_hypertriangulation} in the generic setting.

\begin{remark}
  \label{rmk:general_position}
  According to our current definition of a generic set, $A$, no three points of $\Level{A}{k}$ are collinear, for any $k \in [n-1]$.
  This implies that the projection of any $2$-dimensional face of $\Delta_n^{(k)}$ is a triangle.
  This property also holds if we weaken our notion of genericity to requiring that no three points of $A$ are collinear, and this is the definition we will use from now on.
  However, in this case, two or more points of $\Level{A}{k}$ may coincide, but since they have different labels, we still treat them as different points.
  In the presence of coinciding points, we require that at most one of these points is used in any one triangulation.
  Equivalently, we require that the hypertriangulation remains a hypertriangulation if we perturb the points in $A$ by any sufficiently small amount.
  For an example see Figure~\ref{fig:equalpoints}, which shows two geometrically identical projections of four $2$-dimensional faces of the octahedron, $\Delta_4^{(2)}$.
  Since the respective central vertices have different labels ($13$ versus $24$), the two hypertriangulations are considered different.
  \begin{figure}[htb]
    \centering \vspace{0.1in}
    \resizebox{!}{1.6in}{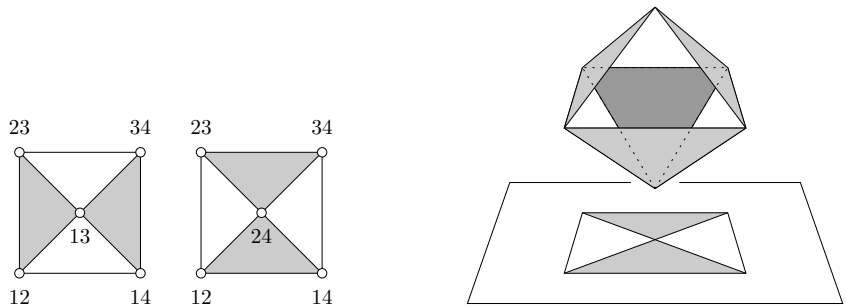}
    \caption{\emph{Left:} two hypertriangulations of four points arranged as the vertices of a diamond in the plane.
    \emph{Right:} one of the hypertriangulations as the projection of an octahedron.}
    \label{fig:equalpoints}
  \end{figure}
\end{remark}

\begin{definition}\label{def:level-k_hypersimplicial_subdivision}
  A \emph{level-$k$ hypersimplicial subdivision} of $A$ is a tiling of $\conv{\Level{A}{k}}$ with projected faces of $\Delta_n^{(k)}$ that remains a tiling under any sufficiently small perturbation of $A$.
\end{definition}

\begin{figure}[htb]
  \centering \vspace{0.1in}
  \resizebox{!}{2.6in}{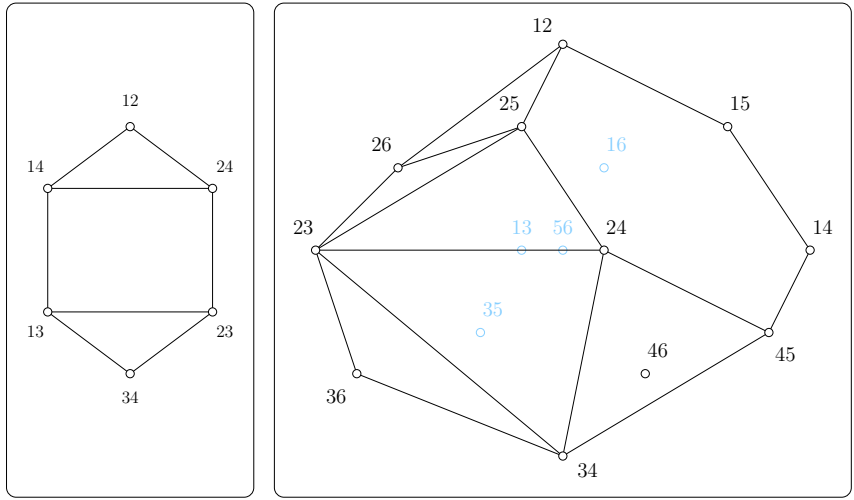}
  \caption{\emph{Left:} a subdivision of the $2$-fold sums of four points that is not a hypersimplicial subdivision of the six points.
  \emph{Right:} a subdivision of the $2$-fold sums of the six points in Figure~\ref{fig:alltypes} that can be interpreted as hypersimplicial in more than one way.}
  \label{fig:nonhyper}
\end{figure}

Note that Definition~\ref{def:level-k_hypersimplicial_subdivision} allows for projections of faces of dimension larger than $2$, which are convex polygons with possibly more than three edges.
In contrast to hypertriangulations, Conditions (V) and (E) of Definition~\ref{def:level-k_hypertriangulation} do not suffice to distinguish hypersimplicial subdivisions from other tilings of $A$.
Take for example a set $A$ consisting of three points spanning a triangle with a fourth point inside this triangle.
Then $\conv{\Level{A}{2}}$ is a convex hexagon, which we may tile as shown in the left panel of Figure~\ref{fig:nonhyper}.
All edges satisfy Condition (E), but this tiling cannot be obtained as projection of faces of the octahedron $\Delta_4^{(2)}$.

\medskip
Similarly, two different hypersimplicial subdivisions may lead to geometrically identical tilings.
See the right panel in Figure~\ref{fig:nonhyper}, in which we cannot distinguish between the lower right triangle being the projection of a tetrahedron or of its $2$-dimensional face whose projected image (which is the triangle) contains the projected fourth vertex.
These two choices are treated as different subdivisions, with the latter being a refinement of the former.
To disambiguate the situation, we will draw the projection of the fourth vertex only if the triangle is the projection of the tetrahedron but not if it is the projection of a $2$-dimensional face of the tetrahedron.

\subsection{Coherent Hypersimplicial Subdivisions and Fiber Polytopes}
\label{sec:2.3}

Within all hypersimplicial subdivisions of a finite set, $A \subseteq \Rspace^2$, the coherent hypersimplicial subdivisions form a privileged subfamily.
These subdivisions are constructed as lower boundaries of convex polytopes obtained by lifting the points of $A^{(k)}$ according to a linear function on $\Rspace^n$.
Let $h \colon \Rspace^n \to \Rspace$ be such a linear function, and write $h_i = h (e_i)$ for the value at the basis vector $e_i$. 
For every $k$, the value at a vertex of $\Delta_n^{(k)}$ is the sum of the values of the corresponding $k$ vertices of $\Delta_n$.
We therefore write $h_I = h (e_I) = \sum\nolimits_{i\in I} h_i$ for every $I \subseteq [n]$.

\begin{definition}
  \label{def:coherent_hypersimplicial_subdivision}
  Let $A_h^{(k)} = \{(a_I,h_I) \mid |I|=k\}$ be the \emph{$h$-lifted} points at level $k$ and note that these are points in $\Rspace^3$.
  The associated \emph{coherent hypersimplicial subdivision}, denoted $T^{(k)}(A,h)$, is the natural projection (which removes the last coordinate) of the lower boundary of $\conv{A_h^{(k)}}$ to $\conv{A^{(k)}}$.
  In the particular case in which the lower boundary has only triangular faces (and no points in the interiors of these faces), we call its projection a \emph{coherent hypertriangulation} of $A$.
\end{definition}

The corresponding \emph{fiber polytope}, denoted $\mathcal F(\Delta_n^{(k)} \to A^{(k)})$, is the set of points
$$
  \frac{1}{\mathrm{area}( \conv{A^{(k)}} )} 
     \int_{ \conv{A^{(k)}} } 
     f(x) \,{\rm d}x,
$$
over all continuous functions $f \colon \conv{A^{(k)}} \to \Delta_n^{(k)}$ 
such that $\pi(f(x))=x$ for all $x\in \conv{A^{(k)}}$.
Note that this is a subset of $\Delta_n^{(k)} \subseteq \Rspace^n$.
The vertices of the fiber polytope correspond to coherent hypertriangulations, and the faces correspond to all coherent subdivisions; see \cite{OlSa21} and \cite{BiSt92} for more details.

\medskip
Because $\Delta_n^{(k)}$ is a complex of hypersimplices, we can use the following equivalent definition, which resembles the one of the secondary polytope associated with the point set $A$; see~\cite{GKZ94}.

\begin{definition}
  \label{def:fiber_polytope}
  Set $e_I = \sum_{i \in I} e_i$ for each $I \subseteq [n]$, and write $\GKZ{\Delta} = \area(\Delta) \cdot \tfrac{1}{3} (e_I+e_J+e_K)$ for every triangle $\Delta = a_Ia_Ja_K$.
  Then the corresponding \emph{fiber polytope}, or \emph{level-$k$ hypersecondary polytope} of $A$,
  denoted $\mathcal F^{(k)}(A)$, is the convex hull of the points
  $$ 
 \GKZ{T} = \frac{1}{\mathrm{area}(  \conv{A^{(k)}} )}  
    \sum\nolimits_{\Delta\in T} \GKZ{\Delta} ,
  $$
  where the points are constructed for all level-$k$ hypertriangulations $T$ of $A$.
\end{definition}

\section{Flips in Hypertriangulations}
\label{sec:3}

The level-$1$ hypertriangulations of a finite set, $A \subseteq \Rspace^2$, are commonly called the (partial) triangulations of $A$.
Flips are the elementary operations that transform one triangulation of $A$ to another.
In the generic situation, there are two types:
the first substitutes one diagonal of a convex quadrangle by the other, and the second subdivides a triangle into three by adding a vertex or coarsens by removing a degree-$3$ vertex; see the Type-I and Type-II flips in Figure~\ref{fig:flips}.
For just four points the flips provide transitions between the only triangulations on these points.
As shown in \cite{HNU99}, these two operations suffice to transform any triangulation to any other triangulation of $A$.

\medskip
This section introduces similar operations for level-$k$ hypertriangulations in $\Rspace^2$.
Before giving the formal definitions, we take a look at configurations of just four points, for which our flips appear naturally.

\subsection{Hypertriangulations of Four Points}
\label{sec:3.1}

\begin{figure}[hbt]
  \centering \vspace{0.1in}
  \resizebox{!}{1.95in}{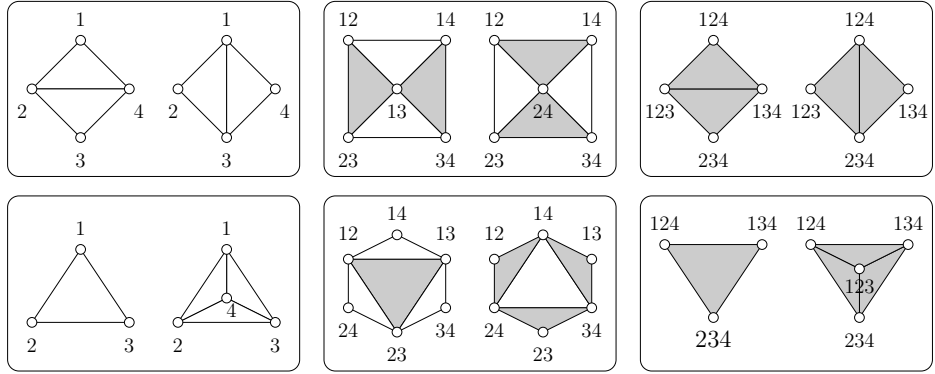}
  \caption{
    The two configurations of four points in $\Rspace^2$ and their hypertriangulations: convex configuration in \emph{top row} and non-convex configuration in \emph{bottom row}.
    \emph{From left to right:} the two level-$1$, level-$2$, and level-$3$ hypertriangulations for each configuration.
    Observe that the squares in the \emph{upper middle} can be more general parallelograms so that the respective central fifth vertices are not at the same geometric location.
    Similarly, the convex hexagons in the \emph{lower middle} are not necessarily regular but are necessarily centrally symmetric.}
  \label{fig:4points}
\end{figure}

For $n = 4$ points in $\Rspace^2$, we have level-$k$ hypertriangulations for $k = 1, 2, 3$.
In the generic case, there are only two combinatorially different configurations of four points:  the vertices of a convex quadrangle, or the vertices of a triangle with the fourth point inside the triangle.
We refer to them as the \emph{convex configuration} and the \emph{non-convex configuration}, respectively.
As illustrated in Figure~\ref{fig:4points}, we have two hypertriangulations for each $k$ and each of the two configurations, so twelve hypertriangulations altogether.
\begin{itemize}
  \item For $k=1$, the vertices are the original points, and all triangles are white.
  \item For $k=2$, there are six points, each the sum of two points in $A$.
  If $A$ is a convex configuration, the convex hull of $\Level{A}{2}$ is a parallelogram and the remaining two points lie inside the parallelogram.
  There are only two hypertriangulations, each choosing one of the two extra points as a vertex and decomposing the parallelogram into two white and two black triangles.
  If $A$ is a non-convex configuration, the points in $\Level{A}{2}$ are the vertices of a centrally symmetric convex hexagon, and there are again only two hypertriangulations.
  \item For $k=3$, the situation is similar to the case $k=1$, except that all triangles are black.
\end{itemize} \medskip
The six pairs of hypertriangulations of four points inspire our definition of flips for hypertriangulations in $\Rspace^2$.

\subsection{Flips}
\label{sec:3.2}

\begin{figure}[hbt]
  \centering \vspace{0.1in}
  \resizebox{!}{1.9in}{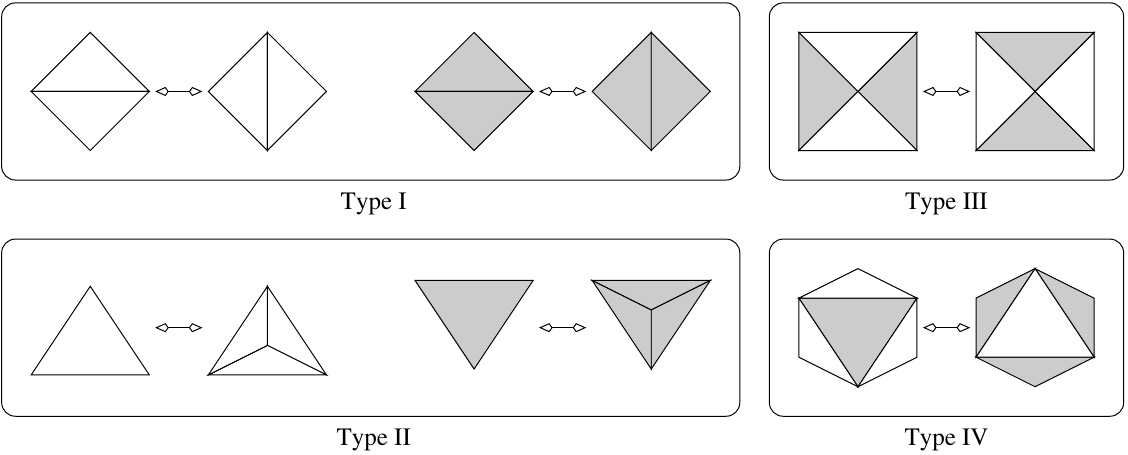}
  \caption{
    The four types of flips in hypertriangulations.
    Type-I and Type-II flips apply to all white or all black triangles, while Type-III and Type-IV flips involve white and black triangles before and also after the flip.}
  \label{fig:flips}
\end{figure}

We introduce four types of flips in hypertriangulations; all illustrated in Figure~\ref{fig:flips}.
A flip preserves the level of the hypertriangulation, so the vertices and edges it introduces must satisfy Conditions (V) and (E) of Definition~\ref{def:level-k_hypertriangulation}.

\begin{definition}
  \label{def:types_I_II_III_IV_flips}
  A \emph{Type-I flip} removes two triangles that share an edge and whose union is a convex quadrangle, and it adds the other two triangles whose union is the same quadrangle.
  The four triangles, two before and two after the flip, are either all white or all black.

  \medskip
  A \emph{Type-II flip} replaces a triangle by three triangles sharing a vertex inside the removed triangle or, in the other direction, replaces the three triangles sharing a degree-$3$ vertex with a single triangle.
  The four triangles, before and after the flip, are either all white or all black.

  \medskip
  A \emph{Type-III flip} applies to a parallelogram decomposed into four triangles, which alternate between white and black around the shared vertex, and replaces these triangles by their reflections through the center of the parallelogram.
  The reflection also switches the colors.

  \medskip
  A \emph{Type-IV flip} applies to a centrally symmetric convex hexagon decomposed into four triangles, one in the middle with a color that is different from the shared color of the surrounding three triangles.
  The flip replaces the four triangles by their reflections through the center of the hexagon, and the reflection switches the colors, as before.
\end{definition}

\begin{remark}
  \label{rmk:redundant_requirements}
  Explicitly requiring colors of triangles and Condition (E) of edges is sometimes excessive.
  For example, if two triangles have the same color, share an edge, and form a convex quadrangle, then Condition (E) is necessarily satisfied by the new edge.
  Similarly, if a vertex is shared by a cycle of four triangles with alternating colors, then these four triangles can be replaced by a Type-III flip.
  On the other hand, for Type-II and Type-IV flips, Condition (E) needs to be taken into account at least in one direction.
\end{remark}

\begin{example}
  \label{ex:all_types_of_flips}
  Take another look at Figure~\ref{fig:alltypes}, which shows a level-$3$ hypertriangulation of six points in $\Rspace^2$.
  It is chosen so that all different types of flips can be applied.
  On the upper left, we see three white triangles that can be replaced by a single white triangle in a Type-II flip.
  Below them, we see three black triangles that can be replaced by a single black triangle in another Type-II flip.
  In the upper middle, we see a hexagonal region with a black triangle surrounded by three white triangles, which can be replaced by a white triangle surrounded by three black triangles in a Type-IV flip.
  After applying this flip, we get a convex quadrangle decomposed into two black triangles, which can be replaced by two other black triangles in a Type-I flip.
  Below that hexagonal region, we see another with a white triangle surrounded by three black triangles, and after applying a Type-IV flip, we get an elongated convex quadrangle decomposed into two white triangles, which can be replace by two other white triangles in another Type-I flip.
  Finally on the right, we see a parallelogram whose triangles alternate between black and white, which can be replaced by four other triangles in a Type-III flip.
\end{example}

\section{Level-2 Hypertriangulations}
\label{sec:4}

This section proves that for every generic point set, $A\subseteq \Rspace^2$, the collection of level-$2$ hypertriangulations is flip-connected.
The main tool used for this purpose is the aging function and its inverse.
This function and its inverse appeared as {\sc UP} and {\sc DOWN} functions for triangulated plabic tilings in the works Galashin \cite{Gal18} and Balitskiy and Wellman \cite{BaWe20} on plabic graphs, and in a modified settings as functions $\mathcal U$ and $\mathcal D$ for zonotopal tilings in the work of Olarte and Santos \cite{OlSa21} on hypersimplicial subdivisions, and as transitions between generations of slices of rhomboid tilings in the work of Edelsbrunner and Osang \cite{EdOs22} on a fast algorithm for level-$k$ Delaunay mosaics.


\subsection{Aging Function}
\label{sec:4.1}

The \emph{aging function}, $F$, maps a white triangle with vertices in $\Level{A}{k}$ to a black triangle with vertices in $\Level{A}{k+1}$.
Specifically, if $\card{I} = \card{J} = \card{K} = k$ and $\card{I \cap J \cap K} = k-1$, then
$$
  F(a_I a_J a_K)  =  a_{I \cup J} a_{J \cup K} a_{K \cup I} .
$$
Indeed, we have $$\card{I \cup J} = \card{J \cup K} = \card{K \cup I} = k+1$$ and $$\card{(I \cup J) \cap (J \cup K) \cap (K \cup I)} = k-1,$$ so the image of $a_I a_J a_K$ under $F$ is black.
The inverse of the aging function is well defined on the black triangles:
if $\card{I} = \card{J} = \card{K} = k$ and $\card{I \cap J \cap K} = k-2$, then
$$
  F^{-1} (a_I a_J a_K) = a_{I \cap J} a_{J \cap K} a_{K \cap I}.
$$
Recall that $W(T)$ and $B(T)$ are the white and black triangles in $T$.
Accordingly, we write $F(W(T))$ and $F^{-1}(B(T))$ for the images under the aging function and its inverse.
With this notation, we have the following property.
\begin{lemma}\label{lem:level_1_to_level_2}
  For every level-$1$ hypertriangulation, $T$, of a generic set $A \subseteq \Rspace^2$, there exists a level-$2$ hypertriangulation, $U$, of $A$ such that $B(U) = F(W(T))$.
\end{lemma}
\begin{proof}
  Since $T$ is level-$1$, all its triangles are white.
  To get $F(W(T))$, we take the midpoints of all edges in $T$ and for every triangle in $T$ draw a (black) triangle that connects the midpoints within this triangle.
  After inflating the configuration by a factor $2$, we get a subset of $\Level{A}{2}$ together with a collection of black triangles; see Figure~\ref{fig:aging}, which suppresses the inflation for better visualization. 
  The part of the convex hull of $\Level{A}{2}$ not covered by black triangles is thus split into regions, and it remains to show that all these regions can be triangulated using white triangles only.
  
  \begin{figure}[hbt]
    \centering \vspace{0.1in}
    \resizebox{!}{1.4in}{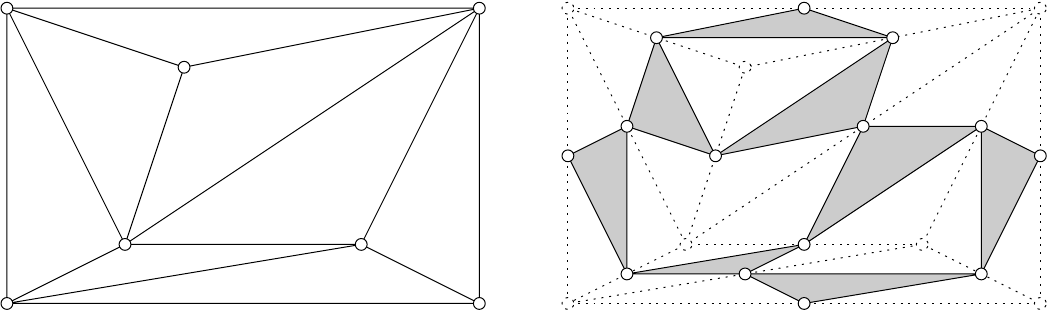}
    \caption{A level-$1$ hypertriangulation of $7$ points on the \emph{left}, and the (shrunken) image of the aging function applied to its $8$ white triangles on the \emph{right}.
    The area left out by the black triangles are (shrunken) copies of the vertex stars in the level-$1$ hypertriangulation.}
    \label{fig:aging}
  \end{figure}

  The gaps between the black triangles are of two types.
  Each gap that is completely surrounded by black triangles is in the shape of the star of an internal vertex of $T$.
  Indeed, for every internal vertex, $a_i$ of $T$, the black triangles within the star of $a_i$ cut out a scaled copy of the said star,
  which inflates into a copy of the star. 
  This star can be triangulated using the vertices in its boundary, and since these vertices share $i$ as one of their labels, all resulting triangles are white.

  \medskip
  Each gap that is not completely surrounded by black triangles is obtained from the star of a boundary vertex of $T$.
  Specifically, for every boundary vertex, $a_j$ of $T$, the black triangles in its star cut out a scaled copy of the star, which, after inflation, is intersected with $\conv{\Level{A}{2}}$.
  This intersection consists of one or more possibly non-convex polygons.
  The vertices of these polygons are the vertices of the star (other than $a_j$), as well as vertices of $\conv{\Level{A}{2}}$.
  The former share $j$ as one of their labels.
  To see that this is also true for the latter, we rotate two parallel supporting lines, one around $\conv{A}$ and the other around $\conv{\Level{A}{2}}$.
  Whenever the second supporting line passes through a boundary vertex that is also a vertex of one of the polygons, the first supporting line passes through $a_j$.
  It follows that $j$ is a label of the boundary vertex.
  Hence, all vertices of the polygons share $j$ as one of their labels, and any triangulation obtained by connecting these vertices produces only white triangles.
\end{proof}

\subsection{Inverse of Aging Function}
\label{sec:4.2}

The aging function can be reversed to construct the level-$1$ hypertriangulation that gives rise to a given level-$2$ hypertriangulation.
To prepare this construction, we let $U$ be a level-$2$ hypertriangulation of $A$, and for each $a_i \in A$, we write $W_i(U)$ for the set of white triangles in $U$, whose three vertices share the label $i$.
In the forward construction of Lemma~\ref{lem:level_1_to_level_2}, $W_i (U)$ would be the triangles re-triangulating (part of) the star of $a_i$, but a priori it is not clear that $U$ can be constructed this way.
We will sometimes abuse notation and write $W_i (U)$ for the union of its triangles.
\begin{lemma}\label{lem:star-convex_white_regions}
  Let $U$ be a level-$2$ hypertriangulation of a generic set $A \subseteq \Rspace^2$, let $a_i \in A$. Assume $W_i (U)$ is non-empty, and let $x$ be a point in its interior.
  Then all triangles in $U$ that have a non-empty intersection with the line segment from $x$ to $2 a_i$ belong to $W_i (U)$, and if $a_i$ lies in the interior of $\conv{A}$, then $W_i(U)$ contains the entire line segment.
\end{lemma}
\begin{proof}
  We prove the case in which $a_i$ lies in the interior of $\conv{A}$.
  The case of a boundary point is easier and omitted.
  To get a contradiction, we assume there is a point $x \in W_i (U)$ such that $W_i (U)$ does not contain the entire line segment from $x$ to $2 a_i$.
  Hence, the line segment crosses the boundary of $W_i (U)$, and we let $a_{ij}$ and $a_{ik}$ be the endpoints of the boundary edge that crosses the line segment closest to $x$.
  Let $L$ be the line that passes through $a_{ij}$ and $a_{ik}$.
  Since $x$ and $2 a_i$ lie on opposite sides of $L$, the points $2 a_j$ and $2 a_k$ lie on the same side of $L$ as $x$.
  It follows that the entire black triangle with vertices $a_{ij}$, $a_{ik}$, $a_{jk}$ lies on this side of $L$.
  But then there are points on the line segment from $x$ to $2 a_i$ outside $W_i (U)$ that are closer to $x$ than the crossing with $L$, which is a contradiction.
  Hence, $W_i(U)$ is star-convex and contains the entire line segment from $x$ to $a_i$.
\end{proof}

We use Lemma~\ref{lem:star-convex_white_regions} to prove that the aging function allows us to go back and forth between level-$1$ and level-$2$ hypertriangulations.
\begin{lemma}\label{lem:level_2_to_level_1}
  Let $U$ be a level-$2$ hypertriangulation of a generic set $A \subseteq \Rspace^2$.
  Then $F^{-1}(B(U))$ is a well defined level-$1$ hypertriangulation of $A$.
\end{lemma}
\begin{proof}
  First we prove that if $a_i$ and $a_j$ are the endpoints of a side of $\conv{A}$, then there is exactly one triangle in $B(U)$ with vertex $a_{ij}$.
  Observe that $a_{ij}$ is necessarily a vertex of $\conv{\Level{A}{2}}$, let $a_{ik}$ and $a_{j\ell}$ be the neighboring boundary vertices, and note that $U$ contains the two edges that connect $a_{ij}$ to these vertices. This can be seen by rotating the supporting line for $\conv{\Level{A}{2}}$ parallel to $a_ia_j$ at $a_{ij}$ around $a_{ij}$; rotation in one direction will give the vertex $a_{ik}$ and rotation in the opposite direction will give the vertex $a_{j\ell}$.
  Traverse the sequence of triangles $a_{ij} u v$ incident to $a_{ij}$ and note that there is necessarily an edge, $uv$, such that $u$ and $v$ neither share $i$ nor $j$ as a label.
  Hence, $a_{ij}uv$ is a black triangle incident to $a_{ij}$.
  If there are two such black triangles, then there is a white triangle between them as two black triangles cannot share a side since vertices of one side of any black triangle in $U$ determine the third vertex of the triangle.
  This white triangle is incident to $a_{ij}$ but neither belongs to $W_i(U)$ nor to $W_j(U)$, which is not possible.

  \medskip
  Next let $\Delta = a_{ij} a_{ik} a_{jk}$ be any black triangle in $U$, and suppose that the edge from $a_i$ to $a_j$ is not a side of $\conv{A}$.
  We prove that there is exactly one other black triangle, $\Delta'$, in $U$ that shares $a_{ij}$, and that $\Delta$ and $\Delta'$ lie on opposite sides of the line that passes through $2 a_i$ and $2 a_j$.
  We consider the case in which both $a_i$ and $a_j$ lie in the interior of $\conv{A}$.
  In all other cases, the argument is similar and omitted.
  In the assumed case, $W_i(U)$ and $W_j(U)$ are both star-convex and meet at $a_{ij}$.
  Traversing the triangles $a_{ij}uv$ incident to $a_{ij}$---starting at $u = a_{ik}$ and ending at $v = a_{j\ell}$ while avoiding $\Delta$---we observe that there must be a second black triangle, $\Delta'$.
  Furthermore, there cannot be three black triangles because $a_{ij}$ has only two labels and can therefore not belong to three white regions.
  The property that $\Delta$ and $\Delta'$ lie on opposite sides of the line passing through $2 a_i$ and $2 a_j$ follows from the star-convexity of $W_i(U)$ and $W_j(U)$ and the fact that these two regions touch at $a_{ij}$, which lies on this line and between $a_i$ and $a_j$.

  \medskip
  When we construct $F^{-1} (B(U))$, we get one white triangle next to every boundary edge of $\conv{A}$ and two non-overlapping white triangles on opposite sides of every non-boundary edge.
  A point in the interior of $\conv{A}$ and sufficiently close to a boundary edge is covered by exactly one triangle in $F^{-1} (B(U))$.
  To move to any other point in $\conv{A}$, we walk from triangle to triangle, and each time we leave a triangle we enter another.
  This implies that almost all points in $\conv{A}$ are covered by exactly one triangle in $F^{-1} (B(U))$.
  It follows that $F^{-1} (B(U))$ is a level-$1$ hypertriangulation of $A$.
  It is unique because the construction via the inverse of the aging function is deterministic.
\end{proof}

\subsection{Flip-connectivity}
\label{sec:4.3}

We are now ready to prove the main result of this section.
Given a generic set $A \subseteq \Rspace^2$, consider the graph whose nodes are the level-$k$ hypertriangulations of $A$, with an arc connecting two nodes if there is a flip that transforms one hypertriangulation to the other.
We call the level-$k$ hypertriangulations \emph{flip-connected} if this graph is connected.
\begin{theorem}\label{thm:level_2_flip-connectivity}
  For every generic point set $A \subseteq \Rspace^2$, the level-$2$ hypertriangulations of $A$ are flip-connected.
\end{theorem}
\begin{proof}
  Let $U$ and $U'$ be two level-$2$ hypertriangulations of $A$, and let $T = F^{-1} (B(U))$ and $T' = F^{-1} (B(U'))$ be the corresponding level-$1$ hypertriangulations, which are possibly partial triangulations of $A$.
  If $T = T'$, then $U$ and $U'$ agree on their black triangles.
  Similarly, the regions of white triangles are the same, but they may be differently triangulated.
  Each such region is a convex or non-convex polygon, and it is known that any two triangulations of a polygon are connected by Type-I and Type-II flips; see \cite{HNU99}.
  {\sloppy
  
  }

  \medskip
  So suppose $T \neq T'$.
  It is also well known that possibly partial triangulations of $A$ are connected by Type-I and Type-II flips; see Figure~\ref{fig:flips}.
  We will show that each Type-I flip on level~$1$ corresponds to a Type-III flip on level~$2$, and each Type-II flip on level~$1$ corresponds to a Type-IV flip on level~$2$.
  So we can perform the flips on the two levels in parallel, but note that Type-I and Type-II flips on level~$2$ are sometimes necessary to enable the Type-III and Type-IV flips.
  
  \begin{figure}[hbt]
    \centering \vspace{0.0in}
    \resizebox{!}{1.1in}{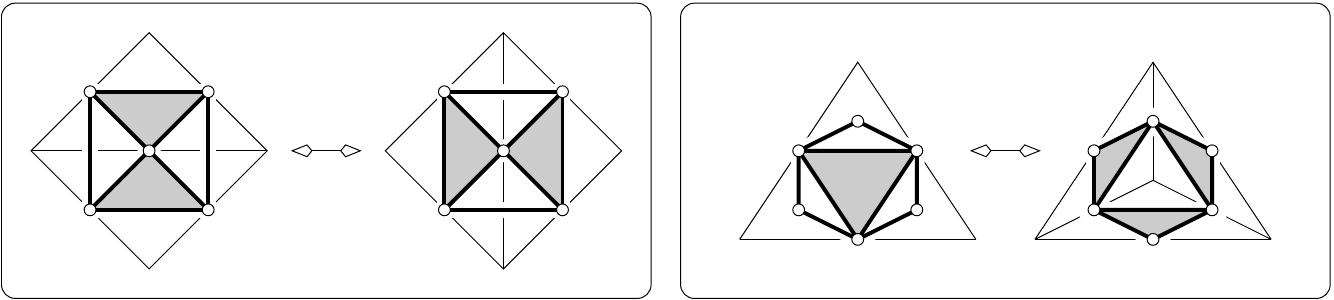}
    \caption{\emph{Left:} a Type-III flip on level~$2$ superimposed on the corresponding Type-I flip on level~$1$.
    \emph{Right:} a Type-IV flip on level~$2$ superimposed on the corresponding Type-II flip on level~$1$.}
    \label{fig:flip-connected}
  \end{figure}
  In the case of a Type-I flip on level~$1$, we need to retriangulate the white regions on level~$2$ as in Figure~\ref{fig:flip-connected} on the left, so a Type-III flip can be performed.
  Such a retriangulation with Type-I flips on level $2$ is always possible.
  The case of a Type-II flip on level~$1$ is similar, except that we need Type-I as well as Tyle-II flips on level $2$ to retriangulate the three white regions around the central black triangle to enable the Type-IV flip; see Figure~\ref{fig:flip-connected} on the right.
  Again, such a retriangulation is always possible.
\end{proof}

We observe that mapping $a_I$ to $a_{[n] \setminus I}$ induces a bijection between the level-$k$ and the level-$(n-k)$ hypertriangulations and their respective flips.
Hence, the flip-connectivity on level $1$ implies the flip-connectivity on level $n-1$.
More interestingly, Theorem~\ref{thm:level_2_flip-connectivity} implies that also the level-$(n-2)$ hypertriangulations of a generic set $A \subseteq \Rspace^2$ are flip-connected.

\subsection{Obstacle for Generalization}
\label{sec:4.4}

\begin{figure}[hbt]
    \centering
    \resizebox{!}{1.6in}{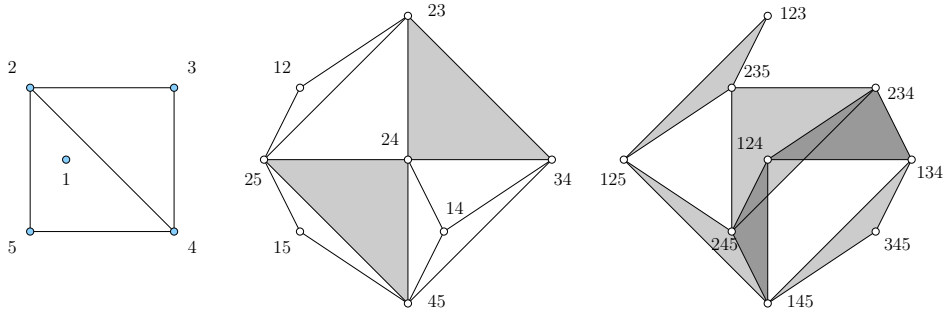}
    \caption{\emph{From left to right:} a partial triangulation of five points with two white triangles, a level-$2$ hypertriangulation with corresponding two black triangles and six additional white triangles, and the corresponding six aged black triangles, some of which overlap.
    The construction is based on Example~5.1 by Olarte and Santos \cite{OlSa21}.}
    \label{fig:noF}
\end{figure}

To summarize, we used the aging function from level $1$ to level $2$ to prove the flip-connectivity of level-$2$ hypertriangulations.
It is not difficult to generalize the aging function to higher levels, but Lemma~\ref{lem:level_1_to_level_2} fails to generalize, which is the reason the authors of this paper were not able to prove flip-connectivity in full generality beyond level $2$.
Indeed, it is known that the extension of Lemma~\ref{lem:level_1_to_level_2} to the aging function that maps white triangles on level $2$ to black triangles on level $3$ fails already for five points.
In particular, the level-$2$ hypertriangulation in \cite[Example~5.1]{OlSa21} contains three triangles, $\Delta_1, \Delta_2, \Delta_3$, whose images under the aging function overlap; see Figure~\ref{fig:noF}.

\subsection{Coherent Hypertriangulations}
\label{sec:4.5}

In this section we briefly discuss the family of coherent level-$k$ hypertriangulations and its flip-connectivity for every $k$.  We note that the aging function and its inverse are well-defined in that setting as for every level-$k$ coherent hypertriangulation $U$ with the height function $h$, the corresponding (non-unique) level-$(k+1)$ and level-$(k-1)$ hypertriangulations with black triangles $F(W(U))$ and white triangles $F^{-1}(B(U))$ may be constructed using the same height function.

\medskip
The flip-connectivity of the family of coherent hypertriangulations relates to properties of the corresponding fiber polytope, $\mathcal F^{(k)}(A)$.
In the simplest case, for $k=1$, the edges of the secondary polytope, $\mathcal F^{(1)}(A)$, correspond to bi-stellar flips, and in the generic case, the flips between coherent triangulations result in edges of the secondary polytope; see \cite[Section~5.3]{DRS10}, where coherent triangulations are called regular.

\medskip
As suggested in \cite[Section~8]{Pos18}, a similar property should hold for all Baues posets (for generic projections) when we use Baues flips instead of bi-stellar flips, but we did not find a precise statement to this effect in the literature.
We therefore include a short sketch of this property for hypertriangulations using the more restrictive family of flips of Types I through IV.
We also refer to \cite[Theorem~5.3.1]{DRS10} for a similar result on usual triangulations (the case $k=1$), which inspired our proof.

\begin{proposition}\label{prop:edge_implies_flip}
  If $\GKZ{T_0}$ and $\GKZ{T_1}$ are two vertices of $\mathcal F^{(k)}(A)$ connected by an edge, then the coherent hypertriangulations $T_0$ and $T_1$ differ by a flip of type I, II, III, or IV (and thus also by a Baues flip).
\end{proposition}

\begin{proof}[Sketch of proof.] 
  Let $T$ be the coherent subdivision that corresponds to the edge connecting $\GKZ{T_0}$ and $\GKZ{T_1}$. 
  Since $T$ is not a hypertriangulation, it contains the projection of a face, $F$, of dimension $3$ or larger.
  If $\dim F\geq 4$, then $F$ is a hypersimplex, and the corresponding coherent hypertriangulations of $F$ give a fiber polytope of dimension at least $2$ and hence cannot be the edge of $\mathcal F^{(k)}(A)$.

  \medskip
  If $T$ contains projections of at least two faces $F$ and $G$ of the hypersimplex $\Delta^{(k)}_n$ with $\dim F= \dim G =3$, then lifting the points of $A^{(k)}$ according to any height function, $h$, from the normal cone of the edge corresponding to $T$ gives two non-triangular faces of the lower boundary of $\conv{A_h^{(k)}}$; here we treat a triangular face with a point in its relative interior as non-triangular.
  
  \medskip
  The labels of vertices $e_I$ of $F$ are obtained by taking 1-, 2-, or 3-element subsets of $\{i,j,k,\ell\}\subset [n]$ and complementing them by a subset of $[n]$ of size $k-1$, $k-2$, or $k-3$. In any case, the lifted points $(a_I,h_I)$ must lie on one two-dimensional plane in $\Rspace^3$ and this results into a linear equation for the heights $h_i,h_j,h_k$, and $h_\ell$. Moreover, the corresponding facet of $\conv{A^{(k)}_h}$ is parallel to the two-dimensional plane through the lifted points $a_i,a_j,a_k$, and $a_\ell$. 
  
  \medskip 
  Similarly, we get a linear equation for the four heights from the face $G$. Since the edge between $\GKZ{T_0}$ and $\GKZ{T_1}$ has codimension $1$, the equations for $F$ and $G$ must be proportional. 
  In that case the corresponding facets of the convex hull of the lifted points are parallel and this is impossible for the lower boundary of $\conv{A_h^{(k)}}$.

  \medskip
  To summarize: the projection of $F$ is the only polygon where $T_0$ differs from $T_1$, and since $\dim F=3$, this results in a flip from $T_0$ to $T_1$.
\end{proof}

Since all coherent level-$k$ hypertriangulations correspond to vertices of the corresponding fiber polytope, the following corollary is immediate.

\begin{corollary}
  \label{cor:coherent_flip-connectivity}
  For every point set $A$, the family of all coherent level-$k$ hypertriangulations is flip-connected.
\end{corollary}

\section{Concluding Remarks}
\label{sec:6}

This section mentions avenues for further research on hypertriangulations and their flips.
In dimension $d = 2$, there is of course the open question of flip-connecti\-vity for levels $k$ between $3$ and $n-3$, in which $n$ is the number of points.

\medskip \noindent \textbf{Beyond 2 dimensions.}
In dimension $d \geq 3$, the question of flip-connectivity for hypertriangulations has yet to be formalized.
Even for generic sets of $n$ points, hypertriangulations beyond level $1$ are generally not simplicial because generic hypersimplices are not necessarily simplices.
The aging function can still be defined and goes through $d$ generations of a $d$-simplex: for $1 \leq j \leq d-1$ from the convex hull of the barycenters of the $(j-1)$-faces to the convex hull of the barycenters of the $j$-faces.
For example in dimension $d=3$, it goes from a tetrahedron (convex hull of its vertices) to an octahedron (convex hull of the midpoints of its edges) to another tetrahedron (convex hull of the barycenters of its triangles).
Flips would be defined as in Section~\ref{sec:3}, which is best described in terms of projections from $d+1$ dimensions.
According to Radon's theorem, there are $d+1$ combinatorially different projections of a $(d+1)$-simplex to $\Rspace^d$ \cite{Rad21}.
The types are paired up, giving $\ceiling{(d+1)/2}$ flips.
The $(d+1)$-simplex has $d+1$ generations, but there is a symmetry between the barycenters of the $(j-1)$-faces and the $(d+1-j)$-faces,
giving $\ceiling{(d+1)/2}$ cases.
We thus get $\ceiling{(d+1)/2}^2$ flips, namely four in $\Rspace^2$, also only four in $\Rspace^3$, but nine in $\Rspace^4$.

%

\medskip \noindent \textbf{Generalized Baues problem.}
One of the main questions in the theory of induced subdivisions is the \emph{generalized Baues problem}.
Roughly, the question asks how well the family of all induced subdivisions embeds the family of coherent subdivisions.
A more specific question asks whether the order complex of the poset of all induced subdivisions is homotopy equivalent to the order complex of the poset of coherent subdivisions.
We refer to the survey of Reiner \cite{Rei99} for precise definitions and details.

\medskip
In the setting of ($2$-dimensional) hypertriangulations, the generalized Baues problem has a positive answer for $k=1$, the case of usual triangulations, as shown by Edelman and Reiner \cite{EdRe98}. Consequently, the problem has also a positive answer for $k=n-1$.
In addition, Olarte and Santos \cite{OlSa21} proved that the generalized Baues problem has a positive answer if the points are in convex position.
For an arbitrary generic set in $\Rspace^2$, the generalized Baues question for level-$k$ hypertriangulations remains open for all $2 \leq k \leq n-2$.

\subsubsection*{Acknowledgements}

Work by all authors but the second is supported by the European Research Council (ERC), grant no.\ 788183, by the Wittgenstein Prize, Austrian Science Fund (FWF), grant no.\ Z 342-N31, and by the DFG Collaborative Research Center TRR 109, Austrian Science Fund (FWF), grant no.\ I 02979-N35.
  Work by the second author is partially supported by the Alexander von Humboldt Foundation and by the Simons Foundation.

\medskip
The second author thanks Jes\'us A. De Loera for useful discussions on flips and non-flips and Pavel Galashin and Alexey Balitskiy for useful discussions on plabic graphs.


\end{document}